\documentclass[journal,twoside,web]{ieeecolor}
\usepackage{amsmath,amssymb,amsfonts}
\usepackage{algorithm}
\usepackage[noend]{algpseudocode}
\usepackage{graphicx}
\usepackage{subcaption}
\usepackage{xcolor}
\usepackage{cite}
\usepackage{bm}
\usepackage{hyperref}
\usepackage{pgfplots}
\pgfplotsset{compat=1.18} 
\usepackage{longtable,tabularx}
\usepackage{wrapfig}
\usepackage{epsfig}
\usepackage{textcomp}
\usepackage{siunitx}
\usepackage{mathrsfs}
\let\labelindent\relax
\usepackage{enumitem}
\usepackage{relsize} 
\usepackage{array}
\definecolor{subsectioncolor}{RGB}{0, 0, 128} 

\newtheorem{theorem}{Theorem}
\newtheorem{remark}{Remark}

\newtheorem{proposition}{Proposition}
\newtheorem{corollary}{Corollary}
\newcommand{\qed}{\hfill$\square$}
\usepackage{fancyhdr}   
\begin{document}

\title{\LARGE \bf
Extremum Seeking for Controlled Vibrational Stabilization of Mechanical Systems: A Variation-of-Constant Averaging Approach Inspired by Flapping Insects Mechanics 
}

\author{{Ahmed A. Elgohary$^{1}$ and Sameh A. Eisa$^{2}$, \IEEEmembership{IEEE Member}}
\thanks{$^{1}$ Department of Aerospace Engineering and Engineering Mechanics, 
        University of Cincinnati, OH, USA
        {\tt\small elgohaam@mail.uc.edu}}%
\thanks{$^{2}$ Department of Aerospace Engineering and Engineering Mechanics, 
        University of Cincinnati, OH, USA
        {\tt\small eisash@ucmail.uc.edu}}%
}

\maketitle
\thispagestyle{empty}

\noindent\textit{This paper has been accepted and published in IEEE Control Systems Letters. DOI: \href{https://doi.org/10.1109/LCSYS.2025.3570767}{10.1109/LCSYS.2025.3570767}}
\vspace{1em}

\pagestyle{fancy}
\fancyhf{}
\lfoot{\scriptsize © 2025 IEEE. Published in IEEE Control Systems Letters. DOI: \href{https://doi.org/10.1109/LCSYS.2025.3570767}{10.1109/LCSYS.2025.3570767}}
\rfoot{\thepage}

\begin{abstract}
This paper presents a novel extremum seeking control (ESC) approach for the vibrational stabilization of a class of mechanical systems \textcolor{black}{(e.g., systems characterized by equations of motion resulting from Newton's second law or Euler-Lagrange mechanics)}.
Inspired by flapping insects mechanics, the proposed ESC approach is operable by only one perturbation signal \textcolor{black}{and can admit generalized forces that are quadratic in velocities.}  
We test our ESC, \textcolor{black}{and compare it against approaches from literature}, on some classical mechanical systems (e.g., mass-spring and an inverted pendulum systems). We also provide a novel, first-of-its-kind, application of the introduced ESC by achieving a 1D model-free source-seeking of a flapping system.    
\end{abstract}


\section{Introduction}
\subsection{Extremum Seeking Feedback Control Systems}
Extremum seeking control (ESC) are model-free, adaptive control techniques that aim at stabilizing a dynamical system about the optimum state of an objective function \cite{ariyur2003real,KRSTICMain,guay2015time}. The main idea of ESC is to apply \textcolor{black}{perturbation signals} to the system. Then, via feedback law(s) involving measurements of the objective function, the ESC loop steers the system towards the extremum; see \cite{scheinker2024100} for a comprehensive review. Classic ESC techniques \cite{KRSTICMain,yilmaz2023exponential} rely on classical averaging methods \cite{khalil2002nonlinear,Maggia2020higherOrderAvg} in analysis and design, requiring the amplitude of the perturbation signals to be small. On the other hand, control-affine ESC systems \cite{DURR2013,VectorFieldGRUSHKOVSKAYA2018,pokhrel2023higher} rely on Lie bracket approximations. In \cite{pokhrel2023higher}, it was shown that said approximations are forms of averaging that require both the amplitude and the frequency of the perturbation signals to be large.  
 
\subsection{Vibrational Stabilization}
\textcolor{black}{Vibrational stabilization is a phenomenon associated with many physical and engineering systems \cite{blekhman2000vibrational}. When done by design, it is considered a powerful control technique that uses high-amplitude, high-frequency perturbation signals to stabilize a system about, usually, an unstable equilibrium point \cite{bullo2002averaging}. Unlike most traditional control techniques, vibrational stabilization can be used in an open-loop setting (i.e., without feedback)} as exemplified in the classical Kapitza inverted pendulum \cite[Section 4.1]{Maggia2020higherOrderAvg} \textcolor{black}{and recent emerging applications (e.g., complex network systems \cite{nobili2023vibrational} and hovering of flapping insects and flapping-wing robots \cite{taha2020vibrational})}. However, \textcolor{black}{due to the fact that open-loop vibrational control may not stabilize the system about a desired equilibrium or optimal condition,} modified approaches (e.g., \cite{tahmasian2018averaging,tahmasian2018averaging1}) have developed feedback frameworks to \textcolor{black}{achieve controlled viabrational stabilization.} All these works \cite{bullo2002averaging,tahmasian2018averaging,tahmasian2018averaging1} utilize for analysis and design a variation-of-constant (VOC) \cite{agravcev1979exponential,bullo2002averaging,Maggia2020higherOrderAvg} averaging approach that is more accustomed to systems with mechanical structure. 
\subsection{\textcolor{black}{Extremum Seeking and Vibrational Stabilization of Mechanical Systems}}
\textcolor{black}{ESC has been applied to a special class of mechanical systems, which only admits generalized forces that are linear in velocities, as follows \cite{suttner2022extremum}: 
\begin{equation}
\ddot{\mathbf{q}} = \underbrace{\mathbf{Y_0}(\mathbf{q}) - \mathbf{R}(\mathbf{q}) \odot \dot{\mathbf{q}}}_{\mathbf{f}(\mathbf{q},\dot{\mathbf{q}})} + \sum_{i=1}^{n} u^i \mathbf{Y_i}(\mathbf{q}),
\label{eq:litrature_mechanical_system}
\end{equation}
where} \( \mathbf{q}, \dot{\mathbf{q}} \in \mathbb{R}^n \) represent the generalized coordinates and velocities, \( \mathbf{f} \in \mathbb{R}^n \) is a vector function representing generalized forces, \(\odot\) denotes the element-wise product operation in a matrix/vector form, \textcolor{black}{$u^i$ and $\mathbf{Y_i}$ represent oscillatory control inputs and their vector fields. The perturbation signals $u^i$ are taken as high amplitude, high frequency sinusoidal functions under certain orthonormality conditions. The aim in \cite{suttner2022extremum} (and its non dissipative extension in \cite{suttner2023extremum}) is to stabilize \eqref{eq:litrature_mechanical_system} about the minimum of an objective function  
$J(\mathbf{q})$ using only its measurements. The authors in \cite{grushkovskaya2021extremum} proposed an ESC framework that uses two high amplitude, high frequency perturbation signals with particular formulations under singular perturbation and quasi-steady state conditions. They showed that their results can be relevant to some mechanical systems such as mass-spring. Additionally, an observer-based ESC approach \cite{observerESC} has been applied to second-order systems, which can be potentially applicable to some mechanical systems. On the other hand, in vibrational stabilization literature, broader classes of mechanical systems have been analyzed where the generalized forces in \eqref{eq:litrature_mechanical_system} are \textit{quadratic} in velocities, enabling vibrational control of systems involving more complex forces, such as but not limited to aerodynamic induced forces. Moreover, one high amplitude, high frequency perturbation signal has been used \textit{only} in the context of open-loop vibrational stabilization \cite{bullo2002averaging,Maggia2020higherOrderAvg}, while two high amplitude, high frequency perturbation signals under restricted phase shift have been used in feedback-based setting \cite{tahmasian2018averaging,tahmasian2018averaging1}.  
}

\subsection{Motivation and Contribution}
\begin{figure}[t]
    \centering
    \includegraphics[width=0.8\linewidth]{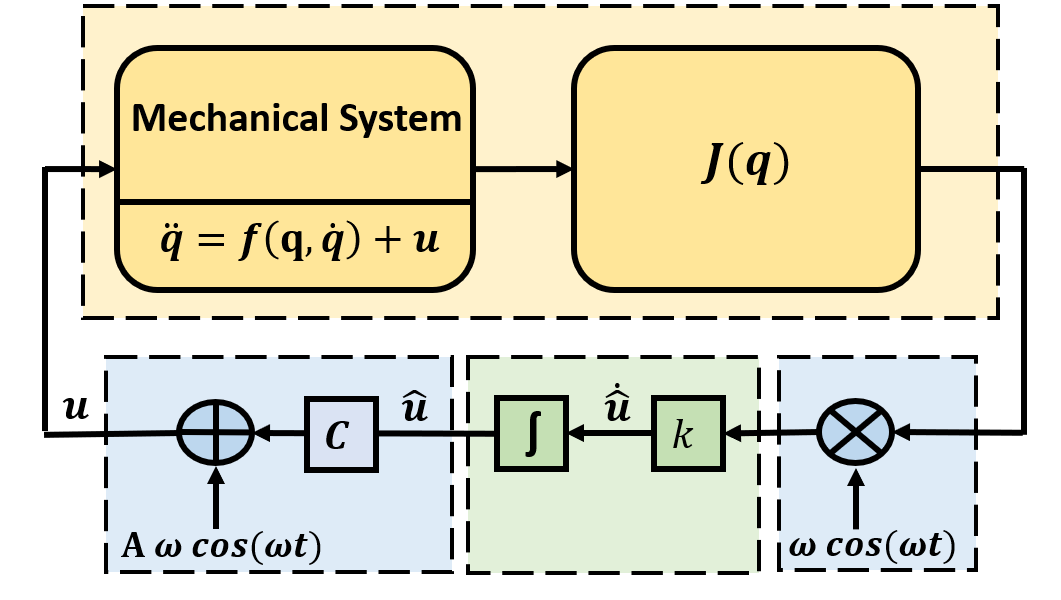}
    \caption{The proposed ESC for a class of mechanical systems.}
\label{fig:esc_framework}
\end{figure}

\textcolor{black}{There are too many mechanical systems in physics, biology, and engineering that can admit (or produce) \textit{only} one source of oscillation/vibration \cite{blekhman2000vibrational}. A natural example of this is that of flapping insects/hummingbirds \cite{taha2020vibrational,enatural_hovering2024}, which enjoy some form of vibrational stabilization using only their wing-flapping action to achieve hovering. The field of experimental biology provides multiple experiments where they also show many flapping insects rely heavily on sensation feedback \cite[Section I.C]{enatural_hovering2024}. This invokes the question: what kind of control method can enable model-free, simple, sensation-based feedback for controlled viabrational stabilization using \textit{only} one perturbation action? Motivated by this question, the authors proposed a special ESC structure (special case of the proposed ESC in this paper in Figure \ref{fig:esc_framework}), which succeeded in characterizing the hovering of flapping insects/hummingbirds using only their wing's natural oscillation, which matched reported experimental data \cite[Table IV]{enatural_hovering2024}. It is important to emphasize, as we clarified in Subsection I.C, that to the best of our knowledge, there is no ESC framework available in literature for mechanical systems that admit generalized forces that are quadratic in velocities (as is the case in aerodynamic forces). Also, none of the available feedback-based ESC or vibrational stabilization methods relevant to mechanical systems are operable with \textit{one} perturbation signal.} 


In this paper, \textcolor{black}{we propose a generalized ESC approach (Figure \ref{fig:esc_framework}) that is inspired by flapping insects mechanics and:
\begin{itemize}
    \item is applicable to many classes of mechanical systems, such as but not limited to systems characterized by equations of motion resulting from Newton's second law and Euler-Lagrange mechanics, with generalized forces allowed to be quadratic in velocities; this is broader than what is available in ESC literature on mechanical systems.
    \item operable by only one perturbation signal, which is: (a) more natural and fits better mechanical systems with one source of vibration (e.g., flapping action in flapping wing flight) -- see also many examples of such systems in \cite{blekhman2000vibrational}; and (b) potentially using less vibration energy as opposed to methods requiring multiple vibration signals.
\end{itemize}
We} analyze the proposed ESC approach using VOC averaging, and show its effectiveness in providing a model-free, real-time control method to viabrationally stabilize a class of mechanical systems about the minimum of an unknown objective function that we have access only to its measurements. 
We test our
approach on some classical mechanical systems \textcolor{black}{with comparison against some literature approaches. We also test successfully the proposed ESC on a system with cubic damping, which suggests the possibility of extending the current approach to even broader classes of mechanical systems.}
Finally, we provide a novel, first-of-its-kind, application of 1D model-free source-seeking of a flapping system.
\section{Main Results}\label{sec:main_results}
Let us consider the following mechanical system:
\small{
\begin{equation}
\ddot{\mathbf{q}} = \mathbf{f}(\mathbf{q}, \dot{\mathbf{q}}) + \boldsymbol{u},
\label{eq:general_mechanical_system}
\end{equation}}
where \( \mathbf{q}, \dot{\mathbf{q}} \in \mathbb{R}^n \) represent the generalized coordinates and velocities, \( \mathbf{f} = [f_1, f_2, \ldots, f_n]^T \in \mathbb{R}^n \) is a vector function that depends on \( \mathbf{q} \) and \( \dot{\mathbf{q}} \), and \( \boldsymbol{u} \in \mathbb{R}^n \) is the control input. We propose one ESC loop (as represented in Figure~\ref{fig:esc_framework}) that is inspired by flapping insects natural hovering ESC mechanics from \cite{enatural_hovering2024}, where there is only one source (e.g., wing flapping oscillation), resulting in only one perturbation signal that influences the system to guide/control its vibrational stabilization about a desired optimal state. The proposed ESC system can be expressed  in first-order state-space representation as:
\small{
\begin{equation}
\frac{d}{dt}
\underbrace{\begin{bmatrix}
\mathbf{q}_{n \times 1} \\
\dot{\mathbf{q}}_{n \times 1} \\
\hat{u}
\end{bmatrix}}_{\textcolor{black}{\mathbf{x}}}
=
\underbrace{\begin{bmatrix}
\dot{\mathbf{q}}_{n \times 1} \\
\mathbf{f}(\mathbf{q}, \dot{\mathbf{q}}) + \mathbf{C} \hat{u} \\
0
 \end{bmatrix}}_{\mathbf{Z}(\mathbf{x})}
+
\underbrace{\begin{bmatrix}
\mathbf{0}_{n \times 1} \\
\mathbf{A} \\
k J(\mathbf{q})
\end{bmatrix} \omega \cos(\omega t)}_{\mathbf{Y}(\mathbf{x}, t)},
\label{eq:proposed_system_dynamics}
\end{equation}}
where \textcolor{black}{\( \mathbf{x} \)} is the state vector that contains $\mathbf{q}$, $\dot{\mathbf{q}}$, and the estimate of control input \( \hat{u} \). The vector fields \( \mathbf{Z}(\mathbf{x}) \) and \( \mathbf{Y}(\mathbf{x}, t) \) represent the drift and oscillatory control effect, respectively. Additionally, \( J(\mathbf{q}) \) represents the unknown objective function to be minimized, which we have access to its measurements and it only depends on \( \mathbf{q} \). The parameter \( k \) is a scaling gain (learning rate) and \( \omega \) is the frequency of oscillation. Additionally, \( \mathbf{C} = [c_1, c_2, \ldots, c_n]^T \in \mathbb{R}^n \) and \( \mathbf{A} = [a_1, a_2, \ldots, a_n]^T \in \mathbb{R}^n \) are vectors of positive gains associated with the estimate of the single control input \( \hat{u} \) and the oscillatory signal \( \omega \cos(\omega t) \), respectively. Now, suppose we have $\mathscr{C}_0$ a closed set and $\mathscr{C}$ an open set, such that $\mathscr{C}_0\subset \mathscr{C} \subseteq \mathbb{R}^{2n+1}$. We assume that the following assumptions hold for the proposed system:
\begin{enumerate}[label=\textbf{A\arabic*}, ref=A\arabic*]
\item The vector fields $\mathbf{Z}(\mathbf{x})$ and $\mathbf{Y}(\mathbf{x}, t)$ and their derivatives up to third-order, are continuous and bounded in $\mathscr{C}_0$.  
    \item The vector function \( \mathbf{f}(\mathbf{q}, \dot{\mathbf{q}}) \) is quadratic in the generalized velocities \( \dot{\mathbf{q}} \), meaning that the third derivative of \( \mathbf{f} \) with respect to \( \dot{\mathbf{q}} \) is zero.
    
    \item The objective function \( J(\mathbf{q}) \) is smooth and has an isolated minimum. That is, there exists a \( \mathbf{q}^* \in \mathscr{C}_0 \) such that \( \nabla J(\mathbf{q}^*) = \bm{0} \), and \( \nabla J(\mathbf{q}) \neq \bm{0} \) for all \( \mathbf{q} \in \mathscr{C}_0 \setminus \{ \mathbf{q}^* \} \). Moreover, with \( J(\mathbf{q}^*) = J^* \in \mathbb{R} \), for all \( \mathbf{q} \in \mathscr{C}_0 \setminus \{ \mathbf{q}^* \} \) there exists two scalar functions $\alpha_1$ and $\alpha_2$ such that $0<\alpha_1(|\bm{q}-\mathbf{q}^*|)\leq J(\mathbf{q})-J^*\leq \alpha_2(|\bm{q}-\mathbf{q}^*|)$.
    \item The vector field \( \mathbf{Y}(\mathbf{x}, t) \) is $T$-periodic in time with zero mean, i.e., \( \int_{0}^{T} \mathbf{Y}(\mathbf{x}, t) \, dt = \bm{0}\), where $T = \frac{2\pi}{\omega}$.
\end{enumerate}
\begin{remark}
    Assumptions A1 and A4 are typical in ESC and vibrational stabilization literature. Assumption A2 means the mechanical system \eqref{eq:proposed_system_dynamics} admits up to quadratic terms in $\dot{\mathbf{q}}$, which is \textcolor{black}{broader than what is available in ESC literature on mechanical systems (e.g., \cite{suttner2022extremum}) and is in level with what is available} in vibrational stabilization of mechanical systems literature (e.g.,  \cite{bullo2002averaging,tahmasian2018averaging,tahmasian2018averaging1}).
    Assumption A3 is common in ESC literature for unknown objective functions with isolated minimum (e.g., \cite{VectorFieldGRUSHKOVSKAYA2018,grushkovskaya2021extremum}) to represent a condition of convexity.   
\end{remark}
\subsection{Variation-of-Constant Averaging Approach}
\begin{figure}
    \centering
    \includegraphics[width=0.6\linewidth]{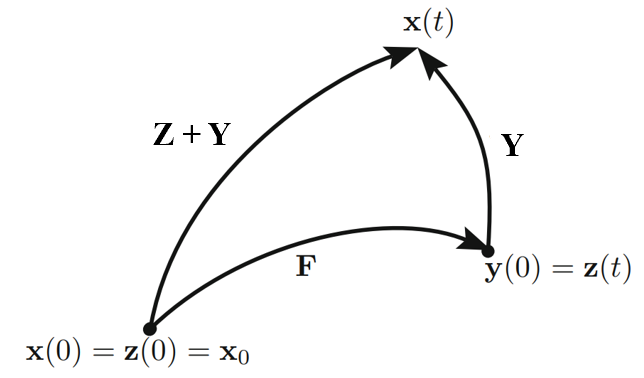}
    \caption{
The flow along \( \mathbf{Z} + \mathbf{Y} \) is equivalent to the flow along the vector fields \( \mathbf{F} \) and then \( \mathbf{Y} \). Starting at \( \mathbf{x}(0) = \mathbf{z}(0) = \mathbf{x}_0 \), the flow first evolves according to \( \mathbf{F} \), and then continues along \( \mathbf{Y} \), resulting in \( \mathbf{y}(0) = \mathbf{z}(t) \).}
    \label{fig:VOC}
\end{figure}

Here we derive the averaged system \textcolor{black}{that captures the qualitative behavior} of our proposed ESC system (\ref{eq:proposed_system_dynamics}) \textcolor{black}{for better characterization and stability analysis}. In fact by considering $\epsilon=1/{\omega}$, (\ref{eq:proposed_system_dynamics}) can be written in the exact form that was shown amenable to VOC-based averaging in \cite{bullo2002averaging,Maggia2020higherOrderAvg}. \textcolor{black}{The reader can refer to \cite{bullo2002averaging,Maggia2020higherOrderAvg} for more details on VOC-based averaging, but we summarize the main idea here. We have two phases. Phase 1, which is concerned with VOC formulation. This phase is independent on averaging, and is based on the idea of using the geometric structure of the control affine system \eqref{eq:proposed_system_dynamics} with two vector fields, i.e., $\mathbf{Z}(\mathbf{x})+\mathbf{Y}(\mathbf{x}, t)$, and find an equivalent representation to it via two separate subsystems along two distinct vector fields, \( \mathbf{F}(\mathbf{z}, t) \) and $\mathbf{Y}(\mathbf{z}, t)$} as depicted in Figure \ref{fig:VOC}; \textcolor{black}{note that $\mathbf{z}$ is a redundant variable used for notational convenience to represent $\mathbf{x}$ in the VOC subsystems.} \textcolor{black}{In phase 2, classical/first-order averaging \cite{khalil2002nonlinear,Maggia2020higherOrderAvg} is applied to the VOC subsystems instead of the original system \eqref{eq:proposed_system_dynamics}. Phase 2 can be particularly useful if one of the VOC subsystems have a zero average as this leads to constant solutions for that subsystem; in fact, this will be the case we have here.} 
The 
two coupled subsystems representing \eqref{eq:proposed_system_dynamics} via VOC \textcolor{black}{(phase 1)} are:
\small{
\begin{equation}
\begin{aligned}
\dot{\mathbf{z}} &= \mathbf{F}(\mathbf{z},t), \quad \mathbf{z}(0) = \mathbf{x}_0, \\
\dot{\mathbf{y}} &= \mathbf{Y}(\mathbf{y}, t), \quad \mathbf{y}(0) = \mathbf{z}(t),
\end{aligned}
\label{eq:subsystems}
\end{equation}}
where \( \mathbf{F} \) represents the pullback of \( \mathbf{Z} \) along the flow of \( \mathbf{Y} \) as first provided in \cite{agravcev1979exponential} (see also \cite{bullo2002averaging,Maggia2020higherOrderAvg}).
\textcolor{black}{For phase 2,} we note here that per assumption A4, \( \mathbf{Y} \) is a zero-mean, time-periodic vector field, resulting in a vanishing average effect \textcolor{black}{given sufficiently small $\epsilon$, or equivalently, large $\omega$ as per first-order averaging \cite{khalil2002nonlinear,Maggia2020higherOrderAvg}}. Consequently, the first-order averaged dynamics of the  system \eqref{eq:proposed_system_dynamics} is determined solely by averaging \( \dot{\mathbf{z}} = \mathbf{F}(\mathbf{z},t)\) \textcolor{black}{in \eqref{eq:subsystems}. For easier representation, we continue our analysis by omitting the dummy variable $\mathbf{z}$ and use $\mathbf{x}$ as the state vector. We proceed by expanding the right-hand side of our VOC system \( \dot{\mathbf{x}} = \mathbf{F}(\mathbf{x},t)\)} as follows \cite{agravcev1979exponential,bullo2002averaging,Maggia2020higherOrderAvg}:
\small{
\begin{multline}
\mathbf{F}(\mathbf{x}, t) = \mathbf{Z}(\mathbf{x}) + \\
\sum_{j=1}^{\infty} \int_0^t \cdots \int_0^{s_{j-1}} \left( \text{ad}_{\mathbf{Y}(\mathbf{x},s_k)} \cdots \text{ad}_{\mathbf{Y}(\mathbf{x},s_1)} \mathbf{Z}(\mathbf{x}) \right) ds_j \cdots ds_1,
\label{eq:pullback_ZY}
\end{multline}}
\textcolor{black}{where $\mathbf{Z}$ and $\mathbf{Y}$ are the vector fields in \eqref{eq:proposed_system_dynamics}} and $\text{ad}_{\mathbf{Y}} \mathbf{Z} = [\mathbf{Y}, \mathbf{Z}]$ represents the Lie bracket, defined 
as:
\small{
\begin{equation}\label{eq:lie_bracket}
[\mathbf{Y}, \mathbf{Z}] = \frac{\partial \mathbf{Z}}{\partial \mathbf{x}} \mathbf{Y} - \frac{\partial \mathbf{Y}}{\partial \mathbf{x}} \mathbf{Z}.
\end{equation}}
\textcolor{black}{For $j=1$, the first term of the series in \eqref{eq:pullback_ZY}, we get $\int_0^t [\mathbf{Y}, \mathbf{Z}](\mathbf{x},s) ds$. For now we compute $[\mathbf{Y}, \mathbf{Z}]$ and leave the integral until a later step.}  
 First, we compute $\frac{\partial \mathbf{Z}}{\partial \mathbf{x}}$ and $\frac{\partial \mathbf{Y}}{\partial \mathbf{x}}$ as:  
\small{
\begin{equation}
\frac{\partial \mathbf{Z}}{\partial \mathbf{x}} =
\begin{bmatrix}
\mathbf{0}_{n \times n} & \mathbf{I}_{n \times n} & \mathbf{0}_{n \times 1} \\
\left[\frac{\partial \mathbf{f}}{\partial \mathbf{q}}\right]_{n \times n} & \left[\frac{\partial \mathbf{f}}{\partial \dot{\mathbf{q}}}\right]_{n \times n} & \mathbf{C}_{n \times 1} \\
\mathbf{0}_{1 \times n} & \mathbf{0}_{1 \times n} & 0
\end{bmatrix},
\end{equation}}
where \( \mathbf{I}_{n \times n} \) is the \( n \times n \) identity matrix and 
\small{
\begin{equation}
\frac{\partial \mathbf{Y}}{\partial \mathbf{x}} =
\begin{bmatrix}
\mathbf{0}_{n \times n} & \mathbf{0}_{n \times n} & \mathbf{0}_{n \times 1} \\
\mathbf{0}_{n \times n} & \mathbf{0}_{n \times n} & \mathbf{0}_{n \times 1} \\
k \nabla_{\mathbf{q}} J(\mathbf{q})_{1 \times n} & \mathbf{0}_{1 \times n} & 0
\end{bmatrix} \omega \cos(\omega t),
\end{equation}
where \( \nabla_{\mathbf{q}} J(\mathbf{q}) \) denotes the gradient of \( J(\mathbf{q}) \) such that \( \nabla_{\mathbf{q}} J(\mathbf{q}) =
\begin{bmatrix} 
\frac{\partial J}{\partial q_1} & \frac{\partial J}{\partial q_2} & \dots & \frac{\partial J}{\partial q_n} 
\end{bmatrix}_{1 \times n} \).} Hence, from \eqref{eq:lie_bracket}:
\small{
\begin{equation}
\label{eq:system9}
[\mathbf{Y}, \mathbf{Z}] =
\begin{bmatrix}
\mathbf{A}_{n \times 1} \\
\left[ \frac{\partial \mathbf{f}}{\partial \dot{\mathbf{q}}} \right]_{n \times n} \left[ \mathbf{A} \right]_{n \times 1} + k J(\mathbf{q}) \left[ \mathbf{C} \right]_{n \times 1} \\
- k \left[ \nabla_{\mathbf{q}} J(\mathbf{q}) \right]_{1 \times n} \left[ \dot{\mathbf{q}} \right]_{n \times 1}
\end{bmatrix} \omega \cos(\omega t).
\end{equation}}
\textcolor{black}{For $j=2$, the second term of the series in \eqref{eq:pullback_ZY}, we get $\int_0^t\int_0^{s_1} [\mathbf{Y},[\mathbf{Y}, \mathbf{Z}]](\mathbf{x},s_1) ds_1ds$. We follow similar procedure to compute $[\mathbf{Y},[\mathbf{Y}, \mathbf{Z}]]= \frac{\partial [\mathbf{Y}, \mathbf{Z}]}{\partial \mathbf{x}} \mathbf{Y} - \frac{\partial \mathbf{Y}}{\partial \mathbf{x}} [\mathbf{Y}, \mathbf{Z}]$ and leave the integral until a later step. In effect, we only need to compute $\frac{\partial [\mathbf{Y}, \mathbf{Z}]}{\partial \mathbf{x}}$ (i.e., taking $\frac{\partial}{\partial \mathbf{x}}$ for \eqref{eq:system9}. Let \(\mathbf{P}_1(\mathbf{q},\dot{\mathbf{q}})=\left[ \frac{\partial \mathbf{f}}{\partial \dot{\mathbf{q}}} \right]_{n \times n} \left[ \mathbf{A} \right]_{n \times 1}\) $+ k J(\mathbf{q}) \left[ \mathbf{C} \right]_{n \times 1}$ and \(\mathbf{P}_2(\mathbf{q},\dot{\mathbf{q}})=-k \left[ \nabla_{\mathbf{q}} J(\mathbf{q}) \right]_{1 \times n} \left[ \dot{\mathbf{q}} \right]_{n \times 1}\) in \eqref{eq:system9}. 
We get the following:}
\small{
\textcolor{black}{
\begin{align}
\frac{\partial [\mathbf{Y}, \mathbf{Z}]}{\partial \mathbf{x}} &= 
\left[
\begin{array}{ccc}
\mathbf{0}_{n \times n} & \mathbf{0}_{n \times n} & \mathbf{0}_{n \times 1} \\[4pt]
\mathbf{M}_{21} &
\mathbf{M}_{22} &
\mathbf{0}_{n \times 1} \\[4pt]
\mathbf{M}_{31} &
- k \nabla_{\mathbf{q}} J(\mathbf{q})_{1 \times n} &
0
\end{array}
\right],
\end{align}}
where \(\mathbf{M}_{21} = \frac{\partial}{\partial \mathbf{q}} (\mathbf{P}_1)\) and \(\mathbf{M}_{31} = \frac{\partial}{\partial \mathbf{q}} (\mathbf{P}_2)\). We note that the terms \(\mathbf{M}_{21}\) and \(\mathbf{M}_{31}\) will be multiplied with zeros in the computation of $\frac{\partial [\mathbf{Y}, \mathbf{Z}]}{\partial \mathbf{x}} \mathbf{Y} - \frac{\partial \mathbf{Y}}{\partial \mathbf{x}} [\mathbf{Y}, \mathbf{Z}]$, so we do not need to compute them in details.
The term \(\mathbf{M}_{22} = \frac{\partial}{\partial \dot{\mathbf{q}}} (\mathbf{P}_1)\) is an $n \times n$ matrix representing second-order mixed derivatives such that for any row index $i=1,...,n$ and column index $k=1,...,n$ of $\mathbf{M}_{22}$, we have the element $\mathbf{M}_{22}|_{(i,k)}=\sum_{j=1}^{n} \frac{\partial^2 f_i}{\partial \dot{q}_k \partial \dot{q}_j} a_j$.
}
\normalsize
\textcolor{black}{Hence, we have:}
\small{
\textcolor{black}{
\begin{equation}
[\mathbf{Y}, [\mathbf{Y}, \mathbf{Z}]] =
\begin{bmatrix}
\mathbf{0}_{n \times 1} \\[6pt]
\left[
\mathbf{M}_{22}\right]_{n \times n}
\left[ \mathbf{A} \right]_{n \times 1} \\[6pt]
- 2k \left[ \nabla_{\mathbf{q}} J(\mathbf{q}) \right]_{1 \times n} \left[ \mathbf{A} \right]_{n \times 1}
\end{bmatrix} \omega^2 \cos^2(\omega t).
\end{equation}}}
\textcolor{black}{For $j=3$}, we follow similar process to the first-order and second-order computations and we find that, 
based on assumption A2, the third-order Lie bracket above and any higher-order Lie bracket (due to the recursive nature of the series \eqref{eq:pullback_ZY}), vanishes. Now, \textcolor{black}{we apply the integrals for $j=1$ and $j=2$. Additionally, we} apply the first-order averaging theorem \cite{khalil2002nonlinear,Maggia2020higherOrderAvg} on \eqref{eq:pullback_ZY}. We get our VOC averaged system: 
\small{
\begin{multline}
\dot{\bar{\mathbf{x}}} = \frac{1}{T} \int_0^T \mathbf{F} dt = \mathbf{Z}(\bar{\mathbf{x}}) + \underbrace{\frac{1}{T} \int_0^T\int_0^t [\mathbf{Y}, \mathbf{Z}](\mathbf{\Bar{x}},s) ds dt}_{I} \\
+ \underbrace{\frac{1}{T} \int_0^T\int_0^t\int_0^{s_1} [\mathbf{Y},[\mathbf{Y}, \mathbf{Z}]](\mathbf{\Bar{x}},s_1) ds_1ds dt}_{II}. 
\end{multline}}
Note that $I$ vanishes due to $cos(\omega t)$ having zero-mean integral with respect to time. The iterated integral of $\omega^2 \cos^2(\omega t)$ in $II$ provide a factor of $1/4$ (see calculation in \cite{MDCLUC2025}). Hence,   
\textcolor{black}{
\begin{equation}
\dot{\bar{\mathbf{x}}} = 
\underbrace{\begin{bmatrix}
\dot{\bar{\mathbf{q}}}_{n \times 1} \\
\mathbf{f}(\bar{\mathbf{q}}, \dot{\bar{\mathbf{q}}}) + \mathbf{C} \hat{\bar{u}} \\
0
\end{bmatrix}}_{\mathbf{Z}(\bar{\mathbf{x}})}
+ \frac{1}{4}
\begin{bmatrix}
\mathbf{0}_{n \times 1} \\[6pt]
\left[
\mathbf{M}_{22}\right]
_{n \times n}
\left[ \mathbf{A} \right]_{n \times 1} \\[6pt]
- 2k \left[ \nabla_{\bar{\mathbf{q}}} J(\bar{\mathbf{q}}) \right]_{1 \times n} \left[ \mathbf{A} \right]_{n \times 1}
\end{bmatrix}.
\label{eq:x_bar_dot}
\end{equation}}


\subsection{Stability Results and Design Guidelines}
\textcolor{black}{As typical in ES methods in literature (e.g., \cite{DURR2013,VectorFieldGRUSHKOVSKAYA2018,pokhrel2023higher}), the stability of our time-varying ESC system \eqref{eq:proposed_system_dynamics} is characterized by its time-invariant averaged system \eqref{eq:x_bar_dot}. As a direct result of the averaging theorem \cite{khalil2002nonlinear,Maggia2020higherOrderAvg}, for small enough $\epsilon$ (or equivalently large enough $\omega$ since $\epsilon=1/{\omega}$), the averaged system \eqref{eq:x_bar_dot} captures the qualitative behavior of \eqref{eq:proposed_system_dynamics} with guarantees of closeness of trajectories as in Corollary \ref{cor:bounded_main}.}
\begin{corollary}\label{cor:bounded_main}
    If $\bar{\mathbf{\textcolor{black}{x}}}(t) \in \mathscr{C}_0$, $\forall t \in [0,t_f/\epsilon]$ with $t_f>0$ and $\bar{\mathbf{\textcolor{black}{x}}}(0)=\mathbf{\textcolor{black}{x}}(0)$, we have $|\mathbf{\textcolor{black}{x}}(t)-\bar{\mathbf{\textcolor{black}{x}}}(t)| = O(\epsilon)$ for $t \in [0,t_f/\epsilon]$.  
\end{corollary}
\begin{remark}
    Corollary 1 provides a bound for the approximation error between the trajectories of $\mathbf{x}$ and the averaged trajectories $\bar{\mathbf{x}}$ in a similar fashion to the approximation errors between the solutions of the original system and its VOC average in vibrational stabilization literature \cite{bullo2002averaging,tahmasian2018averaging,tahmasian2018averaging1}. 
\end{remark}
\textcolor{black}{As a result of Corollary \ref{cor:bounded_main}}, and as guaranteed by the averaging theorem \cite{khalil2002nonlinear,Maggia2020higherOrderAvg}, we have transfer of stability properties from the averaged system \eqref{eq:x_bar_dot} to the original system \eqref{eq:proposed_system_dynamics}; this is also \textcolor{black}{a typical approach} in VOC-based averaging in vibrational stabilization literature \cite{bullo2002averaging,tahmasian2018averaging,tahmasian2018averaging1}. However, to relate the results here more with ESC literature, we provide the following theorem, which guarantees practical stability of \eqref{eq:proposed_system_dynamics} provided that \eqref{eq:x_bar_dot} is asymptotically stable. We refer the reader to \cite[Section 2]{pokhrel2023higher} for the formal definitions of asymptotic stability and practical stability due to space limitations.  
\begin{theorem}\label{thm:rstar}
     If the system \eqref{eq:x_bar_dot} has an equilibrium point $\bar{\mathbf{x}}^*=[\bar{\mathbf{q}}^*,\bm{0},\hat{\bar{u}}^*]^T\in \mathscr{C}_0$ that is asymptotically locally uniformly stable, then the system \eqref{eq:proposed_system_dynamics} is practically uniformly asymptotically stable for $\mathscr{C}_0$.    
\end{theorem}
\textit{Proof:} We assume existence and uniqueness conditions in \cite[Hypothesis 1]{moreau2000practical} are satisfied for the domain $\mathscr{C}_0$. 
From Corollary 1, it follows immediately that \cite[Hypothesis 2]{moreau2000practical} is satisfied for the domain $\mathscr{C}_0$. This proof follows similar rationale to \cite[Theorem 3]{pokhrel2023higher} for $r=1$ (first-order averaging).
\qed
\textcolor{black}{
\begin{remark}
  Theorem 1 means that if one can show the averaged system \eqref{eq:x_bar_dot} is asymptotically stable for some equilibrium point, then the ESC system \eqref{eq:proposed_system_dynamics} is practically stable in the sense that it stays close to that equilibrium point (e.g., oscillates around it). Next, we want to show that \eqref{eq:x_bar_dot} is generally guaranteed to be asymptotically stable, which ensures the ESC system \eqref{eq:proposed_system_dynamics} is generally stable in the sense of Theorem 1.     
\end{remark}}
Given the mechanical structure of the proposed ESC and its amenability to VOC, we expect it may be possible to use Lyapunov function candidates \textcolor{black}{for \eqref{eq:x_bar_dot} to show its asymptotic stability} similar to \cite{bullo2002averaging}. What makes this possible is that by adding an artificial state to the generalized coordinates $\mathbf{q}$ in $\mathbf{x}$, namely $U$ being anti-derivative of $\hat{u}$, we can have a complete mechanical structure similar to those in \cite{bullo2002averaging} with the generalized coordinates being $\mathbf{Q}=[\mathbf{q},U]^T$ and the generalized velocities being $\mathbf{\dot{Q}}=[\mathbf{\dot{q}},\hat{u}]^T$. Nevertheless, we propose a Lyapunov function in the following proposition, taking advantage from assumption A3, which is available in this work as opposed to \cite{bullo2002averaging} due to the ESC framework.
\begin{proposition}
The function $
V = J(\bar{\mathbf{q}}) - J^\star$ is a viable Lyapunov function for \eqref{eq:x_bar_dot}.
\end{proposition} 
\textit{Proof:} From assumption A3, it is clear that $V$ is a positive definite function, which only vanishes at the isolated minimum, i.e., $V(\bar{\mathbf{q}}^*)=0$. Moreover, $\dot{V}=\nabla V \cdot \dot{\bar{\mathbf{x}}}= [\nabla_{\bar{\mathbf{q}}}J,\bm{0}_{1\times n,0}]\cdot[\frac{d\bar{\mathbf{q}}}{dt}^T,\frac{d\dot{\bar{\mathbf{q}}}}{dt}^T,\frac{d\bar{\hat{u}}}{dt}]$. That is, $\dot{V}=\nabla_{\bar{\mathbf{q}}} J \cdot \dot{\bar{\mathbf{q}}}$. From assumption A3, we immediately obtain that \( \dot{V} = 0 \) at \( \bar{\mathbf{q}} = \mathbf{q}^\star \). Now it remains to show that \( \dot{V} < 0 \) for any \( \bar{\mathbf{q}} \neq \mathbf{q}^\star \). We note from \eqref{eq:x_bar_dot} that $\frac{d(\hat{\bar{u}}-\hat{\bar{u}}^*)}{dt}=\frac{d\hat{\bar{u}}}{dt}$ in \eqref{eq:x_bar_dot} has a vector field that is in the gradient descent direction (negative of the gradient). Then, based on assumption A3, $\hat{\bar{u}}-\hat{\bar{u}}^*$ is decreasing (i.e., distance between $\hat{\bar{u}}$ and $\hat{\bar{u}}^*$ is shrinking) near the minimum with steady state $\hat{\bar{u}}=\hat{\bar{u}}^*$ when $\bar{\mathbf{q}}=\mathbf{q}^\star$. Now, we recall that $\dot{\bar{\mathbf{q}}}^*=\bm{0}$ and that
\small{
\begin{equation}
\frac{d(\dot{\bar{\mathbf{q}}}-\dot{\bar{\mathbf{q}}}^*)}{dt} = \frac{d\dot{\bar{\mathbf{q}}}}{dt}=
\underbrace{\bm{f}(\bar{\mathbf{q}}, \dot{\bar{\mathbf{q}}})}_{I_1} 
+ \underbrace{\mathbf{C} {\hat{u}}}_{I_2} 
+ \underbrace{\frac{1}{4}\mathbf{M}_{22} \mathbf{A}}_{I_3}.
\label{eq:qbar_dot}
\end{equation}}
Given the boundedness of \( I_1 \), \( I_3 \), and \( \hat{u} \) from assumption A1, we argue that there exist large entries \( c_1, c_2, \dots, c_n \) such that \( C\hat{u} \) dominates over \( I_1 \) and \( I_3 \). With $\hat{\bar{u}}-\hat{\bar{u}}^*$ decreasing, we have one of two possibilities: (1) \( I_1+I_3 \) is negative near equilibrium (have a stabilizing effect), which means $\hat{\bar{u}}^*$ is positive or zero; or (2) \( I_1+I_3 \) is positive near and at the equilibrium (have a destabilizing effect), which means $\hat{\bar{u}}^*$ is negative.    
In both (1) and (2), due to the dominance of \( C\hat{u} \) over \( I_1+I_3 \), the vector field of $\frac{d(\dot{\bar{\mathbf{q}}}-\dot{\bar{\mathbf{q}}}^*)}{dt}=\frac{d(\dot{\bar{\mathbf{q}}}-\bm{0})}{dt}$ will be negative near equilibrium. Thus, $\dot{\bar{\mathbf{q}}}$ is decreasing in the gradient descent direction (negative of the gradient). Hence, \(  \dot{V}=\nabla_{\bar{\mathbf{q}}} J \cdot \dot{\bar{\mathbf{q}}} < 0 \) for any \( \bar{\mathbf{q}} \neq \mathbf{q}^\star \) in some neighborhood of $\mathbf{q}^\star$.
\qed
\\
\textcolor{black}{\textbf{Design guidelines \& comments.} Since the viability of the averaging theorem is crucial to obtain a stable behavior of the ESC system (i.e., ensuring the viability of Corollary \ref{cor:bounded_main} and Theorem \ref{thm:rstar}), the frequency $\omega=1/\epsilon$ needs to be sufficiently large as a result of $\epsilon$ needs to be sufficiently small. Moreover, $k$ is a tuning parameter that depends on the local behavior of $J(\mathbf{q})$ near the equilibrium as is the case in gradient-based ESC literature; usually the user is recommended to start with smaller values of $k$ when testing the learning rate for $\nabla J(\mathbf{q})$ before resorting to larger values in a similar manner to gradient-based or decent algorithms. Based on Proposition 1, the gains $\mathbf{C}$ and $\mathbf{A}$ need to be balanced in the following way: pick $\mathbf{C}$ to be larger as needed compared to $\mathbf{A}$, which are set to be small as needed (does not have to be very small, but can be). This balance enables the control input of the ESC to overpower/adjust forces represented in the drift (similar observations were made in \cite{suttner2022extremum,suttner2023extremum}). Another important guideline to note, is that the gains in $\mathbf{A}$ (see last equation in \eqref{eq:x_bar_dot}) can be used to leverage some components of $\nabla J(\mathbf{q})$, i.e., make some gradient components larger or smaller; this can be useful for balance if one or more components of $\nabla J(\mathbf{q})$ are very small or very large.}

\section{Applications}
In this section, we demonstrate the effectiveness of the proposed ESC as a model-free method for guided/controlled vibrational stabilization using one perturbation signal. We test our ESC on two classical mechanical systems \textcolor{black}{with comparisons against literature}. Then, we test our ESC in a novel application of 1D source (height) seeking by a flapping system emulating flapping insects/hummingbirds/robots. 
\subsection{Classic Mass-Spring-Damper System}
\begin{figure}
    \centering
    \includegraphics[width=0.85\linewidth]{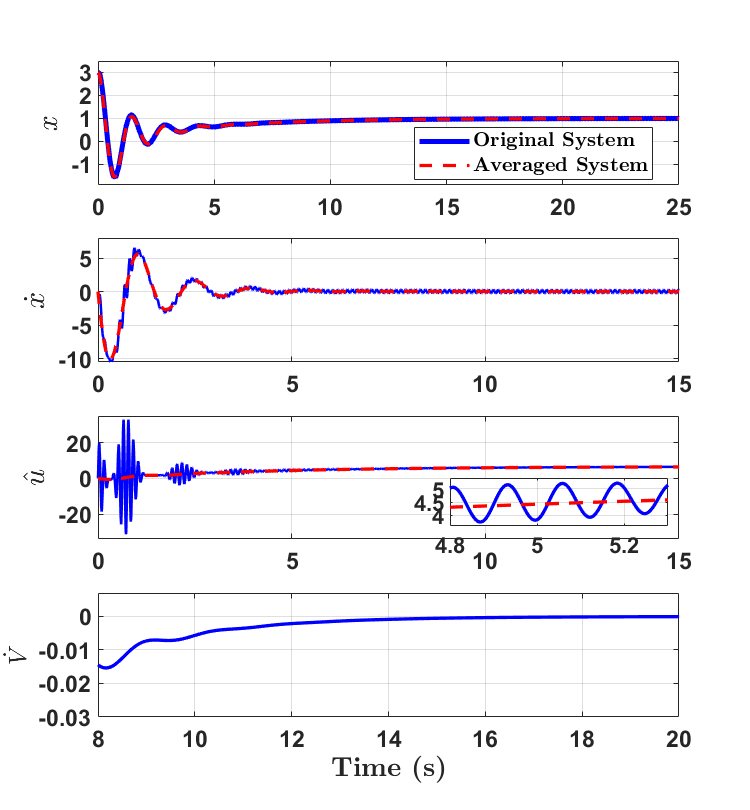}
    \caption{Top 3 subplots: Trajectories of mass-spring system generalized coordinate, its velocity and ESC input, along with their VOC average. Bottom subplot: Lyapunov function
rate.}
    \label{fig:mass_spring}
\end{figure}

We consider the classical mass-spring-damper system, which is known to be stable under its own drift, which represents non-conservative forces. Said system, consisting of a mass \( m \) connected to a fixed support via a linear spring with stiffness \( \alpha \) and a damping element with coefficient \( \beta \), is written in the form \eqref{eq:general_mechanical_system} as follows: 

\begin{equation}\label{eq:mass}
\ddot{x} =
\underbrace{-\frac{\alpha}{m} x - \frac{\beta}{m} \dot{x}}_{f(x, \dot{x})} + u, \quad \text{with} \: u = c_1 \hat{u} + a_1 \omega \cos(\omega t),
\end{equation}
where \( x \) is the displacement of the mass (generalized coordinate), \( \dot{x} \) and \( \ddot{x} \) are its generalized velocity and acceleration. Per the ESC loop (Figure \ref{fig:esc_framework}), the control input \( u \) consists of a high amplitude, high-frequency perturbation signal \( a_1 \omega \cos(\omega t) \) added to the adaptive control estimate \( c_1 \hat{u} \). The adaptation law for the control estimate is:  
\begin{equation}
\dot{\hat{u}} = k J(x) \omega \cos(\omega t),
\end{equation}
where $J(x)$ is the unknown, but measurable, objective function. In our simulation (see Figure \ref{fig:mass_spring}), we take $J(x) = (x - x_{\text{desired}})^2$, where \( x_{\text{desired}} \) represents the optimal equilibrium about which the mass-spring system vibrationally stabilizes. We set \( x_{\text{desired}} = 1 \), \( m = 1 \) kg, \( \alpha = 20 \), and \( \beta = 2 \) Ns/m. The ESC perturbation parameters are chosen as \textcolor{black}{\( c_1 = 3 \), \( a_1 = 0.3 \)}, perturbation frequency \( \omega = 50 \), and adaptation gain \( k = 5 \). The system is initialized as: \( x(0) = 3 \), \( \dot{x}(0) = 0 \), and \( \hat{u}(0) = 0 \). The simulation results, shown in Figure~\ref{fig:mass_spring}, confirm that the proposed ESC successfully achieves vibrational stabilization, guiding \( x \) towards \( x_{\text{desired}} \), adjusting the effect of the stabilizing drift $f(x, \dot{x})$ in \eqref{eq:mass}; this is why $\hat{\Bar{u}}^*$ is positive as observed in the proof of Proposition 1. Calculations for the VOC averaged system \eqref{eq:x_bar_dot} for this example is available in \cite{MDCLUC2025}. 

\begin{figure}
    \centering
    \includegraphics[width=0.8\linewidth]{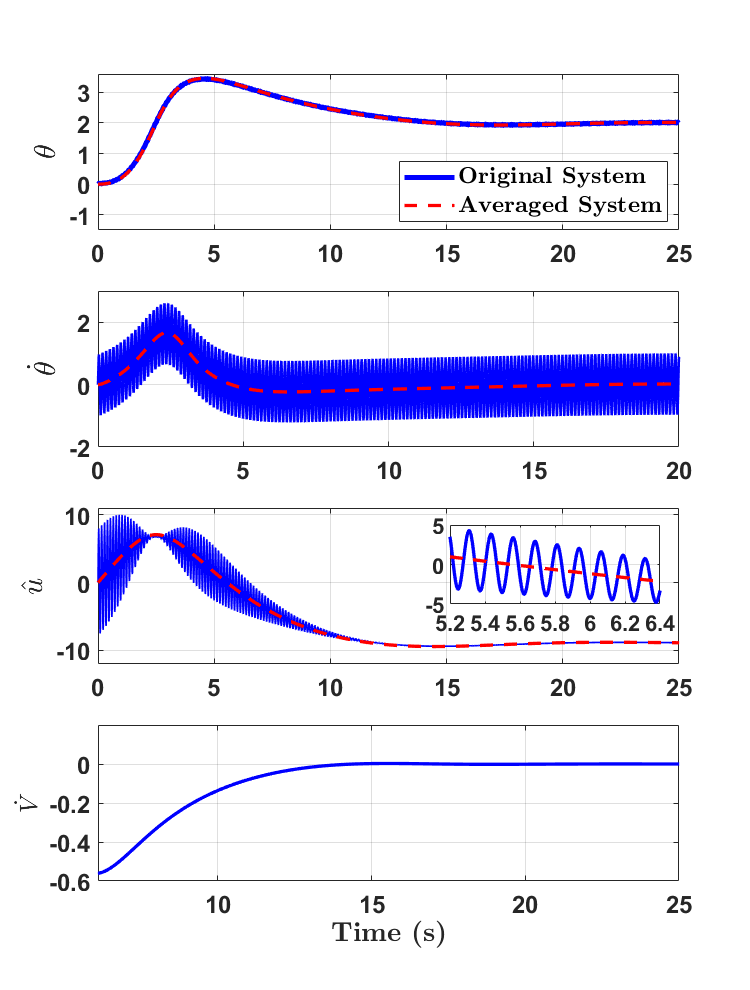}
    \caption{Top 3 subplots: Trajectories of inverted pendulum generalized coordinate, its velocity and ESC input, along with their VOC average. Bottom subplot: Lyapunov function
rate.}
    \label{fig:inverted_pendulum}
\end{figure}

\subsection{Classic Inverted Pendulum}
We consider the classical problem of the inverted pendulum, a known system in the nonlinear control field due its instability under drift, but its ability to be stabilized in open-loop via vibrational stabilization methods \cite{Maggia2020higherOrderAvg}. The system consists of a pendulum of mass \( m \) and length \( L \) pivoting about a fixed point under gravity and damping. Its dynamics are written in the form of (\ref{eq:general_mechanical_system}) as follows: 
\begin{equation}\label{eq:pend}
\ddot{\theta} = 
\underbrace{\frac{g}{L} \sin(\theta) - \frac{\beta}{mL} \dot{\theta}}_{f(\theta, \dot{\theta})} + u, \quad \text{with} \: u = c_1 \hat{u} + a_1 \omega \cos(\omega t),
\end{equation}
where \( \theta \) is the angular position (generalized coordinate), \( \dot{\theta} \) and \( \ddot{\theta} \) are its generalized velocity and acceleration, \( g \) is gravitational acceleration, and \( \beta \) is the damping coefficient. Similar to the mass-spring system (Subsection III.A), we have: 
\begin{equation}
\dot{\hat{u}} = k J(\theta) \omega \cos(\omega t),
\end{equation}
where $J(\theta)$ is the objective function, which we take as $J(\theta) = (\theta - \theta_{\text{desired}})^2$ in our simulation (see Figure \ref{fig:inverted_pendulum}), with \( \theta_{\text{desired}} \) representing the optimal target angle the inverted pendulum is to viabrationally stabilize about. In the simulations, we set \( \theta_{\text{desired}} = 2 \),  \( m = 1 \) kg,  \( L = 1 \) m, \( g = 9.81 \) m/s\(^2\), and  \( \beta = 10 \) Ns/m. The ESC perturbation parameters are chosen as \textcolor{black}{\( c_1 = 1 \), \( a_1 = 0.5 \)}, perturbation frequency \( \omega = 50 \), and adaptation gain \( k = 2 \). The system is initialized as: \( \theta(0) = 0 \), \( \dot{\theta}(0) = 0 \), and \( \hat{u}(0) = 0 \). The simulation results, shown in Figure~\ref{fig:inverted_pendulum}, confirm that the proposed ESC successfully achieves vibrational stabilization, guiding \( \theta \) towards \( \theta_{\text{desired}} \), overcoming destabilizing drift $f(\theta, \dot{\theta})$ in \eqref{eq:pend}; this is why $\hat{\Bar{u}}^*$ is negative as observed in Proposition 1 proof. Calculations for the VOC averaged system \eqref{eq:x_bar_dot} for this example is available in \cite{MDCLUC2025}.

\subsection{\textcolor{black}{Comparison with ESC Literature}}

\textcolor{black}{In this subsection, we try to compare the performance of the proposed ESC with two ESC schemes relevant to mechanical systems: the Lie bracket-based method by \cite{grushkovskaya2021extremum} and the two-dither controller by \cite{suttner2022extremum}. For a fair and meaningful comparison, we maintained the same (or kept our ESC at a lower, if necessary) initial convergence rate. We also used the same perturbation frequency and objective function. We also provide codes for the implementation/results in \cite{MDCLUC2025}.}

\begin{figure}
    \centering
    \includegraphics[width=0.8\linewidth]{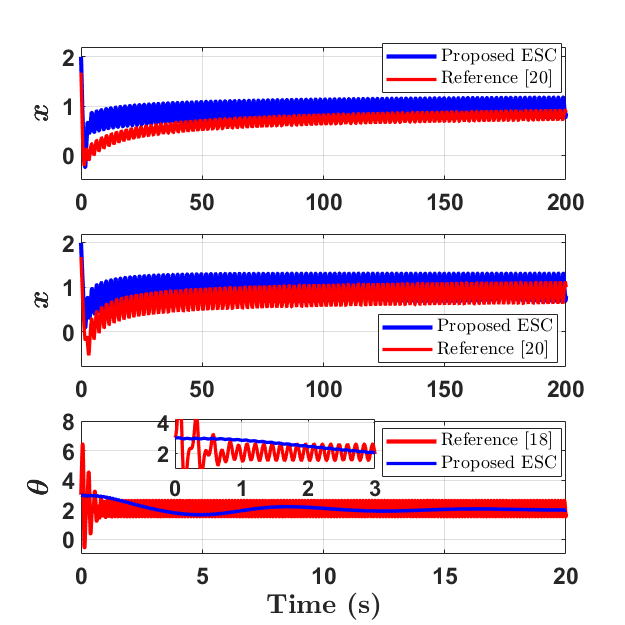}
    \caption{Comparison between the proposed ESC and existing methods from~\cite{grushkovskaya2021extremum} and~\cite{suttner2022extremum}. Top: Mass-spring-damper system with linear damping. Middle: Mass-spring-damper system with cubic damping. Bottom: Inverted pendulum system.}
    \label{fig:comparison}
\end{figure}

\textcolor{black}{\textbf{Comparison with \cite{grushkovskaya2021extremum}.} We implemented the control law in \cite[Eq(6)]{grushkovskaya2021extremum} on the mass-spring system \eqref{eq:mass}. That is, $\dot{u}=\frac{2\sqrt{\pi}}{\eta \sqrt{\epsilon}}\left(g_1(J(x))\cos(\frac{2\pi t}{\eta \epsilon })+g_2(J(x))\sin(\frac{2\pi t}{\eta \epsilon })\right)$ with $g_1(J(x)) = \sqrt{\gamma}J(x)$, $g_2(J(x)) = \sqrt{\gamma}$. To maintain the same initial convergence rate, we conducted the comparative simulations in Figure \ref{fig:comparison} using $J(x) = (x - x_{\text{desired}})^2$, $x_{\text{desired}}=1$, $x(0) = 2$, $\dot{x}(0) = -1$, and $\hat{u}(0) = 0$ for our ESC and using $x(0) = 1.68$, $\dot{x}(0) = -1$, and $u(0) = 0$ for \cite{grushkovskaya2021extremum}. We used parameters $c_1 = 6$, $a_1 = 1$, $\omega = 3.2$, and $k = 0.72$ for our ESC, and used 
$k = 10$, 
$\mu = 5$ (see \cite[Eq(5)]{grushkovskaya2021extremum}) with $\gamma = 100$, $\varepsilon = 1/12.3$, and $\eta = 25$ to keep $\frac{2\pi}{\eta \epsilon }=\omega=3.2$. The comparison was repeated for nonlinear damping by replacing the linear damping term with a cubic one in (\ref{eq:mass}), i.e., replace \( -\frac{\beta}{m} \dot{x} \) with \( -\frac{\beta}{m} \dot{x}^3 \), keeping all other setting unchanged.}

\textcolor{black}{\textbf{Comparison with \cite{suttner2022extremum}.} We implemented the control law in~\cite[Eq(9)]{suttner2022extremum} (see also \eqref{eq:litrature_mechanical_system}) on the inverted pendulum system (\ref{eq:pend}) using $J(\theta) = (\theta - \theta_{\text{desired}})^2$, $\theta_{\text{desired}}=2$. That is, $\mathbf{f}(\mathbf{q},\dot{\mathbf{q}})$ in \eqref{eq:litrature_mechanical_system} is $f(\theta,\dot{\theta})$ in \eqref{eq:pend} with $u_1 = \omega \sin(\omega t)\lambda_1 J(\theta) + \mu_1 \alpha^2(J(\theta))$ and $u_2 = \omega \cos(\omega t)\lambda_2 \alpha(J(\theta)) + \mu_2 \alpha^2(J(\theta))$, $Y_1=Y_2=1$, where $\alpha(J(\theta))=\sqrt{J(\theta)+\log(2\cosh(J(\theta))}$.
We conducted a comparative simulation with our ESC as in Figure \ref{fig:comparison}. Both \eqref{eq:litrature_mechanical_system} and the proposed ESC have $\theta(0) = 3$ and $\dot{\theta}(0) = -1$, and used the same perturbation frequency $\omega = 50$. Since the ESC in \cite{suttner2022extremum} needs two perturbation signals with high amplitudes in the order of $\omega$, we were not able to match the initial convergence rate as the proposed ESC needs much less perturbation effort and operates with much less initial convergence rate (the case here). We also note that we tried to keep the initial convergence rate and parameters of \cite{suttner2022extremum} as small as possible. Our ESC parameters were $c_1 = 7$, $a_1 = 1$, and $k = 1$. For \cite{suttner2022extremum},   
the parameters were $\lambda_1 = \lambda_2 = 14$ and $\mu_1 = \mu_2 = 14$.}

\textcolor{black}{\textbf{Comments.} In both the linear and cubic damping scenarios, our proposed ESC achieved faster convergence than \cite{grushkovskaya2021extremum}, in a relatively lower frequency range similar to what was reported in the simulations in \cite{grushkovskaya2021extremum}; for example, we used much higher frequency on our version of the mass-spring problem in Subsection III.A. The ability of our ESC to handle generalized forces that are cubic in velocities is to be noted as this suggests the potential of extending analysis provided in this paper to ESC admitting more relaxed considitons on the generlized forces. Also, in comparison with \cite{suttner2022extremum}, our method resulted in significantly lower oscillations throughout, including around the equilibrium. Even though our ESC succeeded in solving the inverted pendulum problem with cubic damping, we did not include a comparison with \cite{suttner2022extremum} since that approach only admits generalized forces linear in velocities -- see \eqref{eq:litrature_mechanical_system}. However, we would like to note that there could be systems where \cite{suttner2022extremum} is more suited to solve, especially if the problem involves structured/set $\mathbf{Y}_i$ to be associated with the perturbation signals.}

\subsection{1 Dimensional Source (Height) Seeking by a Flapper}
ESC methods have been introduced successfully to characterize/emulate optimized energy-efficient flight strategies by soaring birds (e.g., \cite{eisa2023}). Following that, we were able to characterize/emulate stable hovering against disturbances in flapping insects/hummingbirds by ESC \cite{enatural_hovering2024}. However, the ESC in \cite{enatural_hovering2024} cannot source/height-seek due to lack of projected force in the vertical acceleration equation of motion. That is, the special case of our ESC (Figure \ref{fig:esc_framework}) in \cite{enatural_hovering2024} can not transition heights to stabilize (hover) about particularly desired height. This not to be the case with the proposed ESC in this paper.  

\begin{figure}
    \centering
    \includegraphics[width=0.85\linewidth]{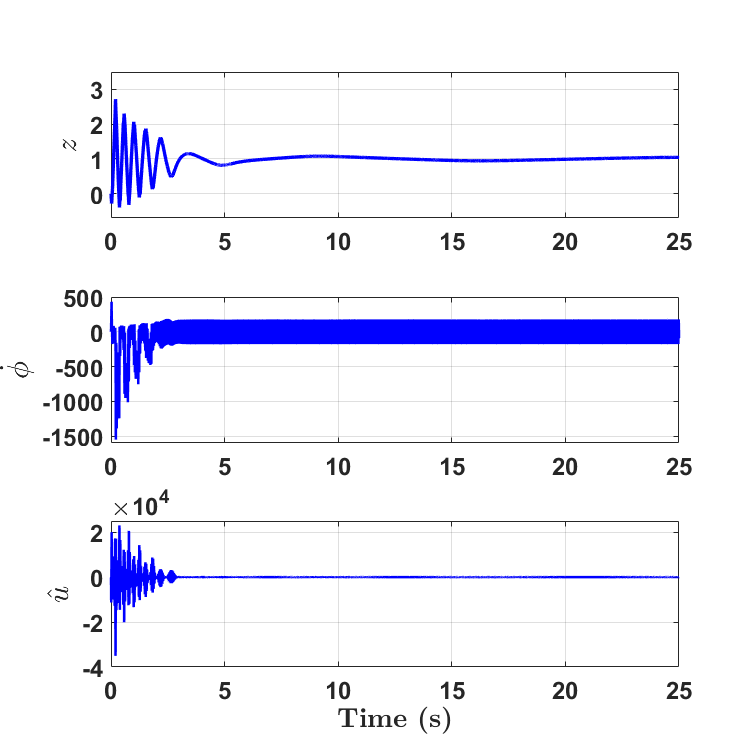}
    \caption{Our ESC successfully stabilizes the flapping system about a desirable height (hover at a certain desirable height). }
    \label{fig:flapping}
\end{figure}

We proceed by adopting 2DOF flapping-wing and aerodynamic model, representing vertical movement, due to multiple advantages that we layout by detail in \cite[Section 2]{enatural_hovering2024}. The considered flapping system is described by the altitude \( z \) and wing flapping angle \( \phi \), with forces and torques acting through aerodynamic effects. This system is written in the form \eqref{eq:general_mechanical_system} as:
\begin{equation}\label{eq:flap}
\begin{bmatrix}
\ddot{z} \\
\ddot{\phi}
\end{bmatrix}
=
\underbrace{
\begin{bmatrix}
- k_{d1} |\dot{\phi}| \dot{z} + g - k_L \dot{\phi}^2 \\
- k_{d3} \dot{z} \dot{\phi} - k_{d2} |\dot{\phi}| \dot{\phi}
\end{bmatrix}
}_{\bm{f}(z, \dot{z}, \phi, \dot{\phi})}
+ 
\mathbf{u},
\end{equation}
where the control input \( \mathbf{u} \) follows the proposed ESC framework as depicted in Figure \ref{fig:esc_framework}:
\begin{equation}
\mathbf{u} =
\begin{bmatrix}
c_1 \\
c_2
\end{bmatrix} \hat{u}
+
\begin{bmatrix}
\frac{a_1}{\textcolor{black}{I_F}} \\
\frac{a_2}{\textcolor{black}{I_F}}
\end{bmatrix} \omega \cos(\omega t).
\end{equation}
\textcolor{black}{where \( I_F \) is the wing
flapping moment of inertia.} The system parameters \cite[Section 2]{enatural_hovering2024} are given by \textcolor{black}{\( I_F = 1.3179 \times 10^{-7} \)}, \( k_{d1} = 0.0353739 \), \( k_L = 0.000621676 \), \( k_{d2} = 0.33915 \), and \( k_{d3} = 16.5766 \). The flapping frequency is set to \( f = 26.3 \) Hz (this is the natural frequency of Hawkmoth insects), with \( \omega = 2\pi f \). The adaptation law for the control estimate is:
\begin{equation}
\dot{\hat{u}} = k J(z) \omega \cos(\omega t),
\end{equation}
where $J(z)$ is the objective function sensed/measured by the insect/mimicking-robot and it depends on height. We take $J(z) = (z - z_{\text{desired}})^2$. In our simulation (see Figure \ref{fig:flapping}), we take \( z_{\text{desired}} = 1 \) to represent the target/source at an altitude. The ESC parameters are selected as \( c_1 = 0.15 \), \( c_2 = 1 \), \textcolor{black}{\( a_1 = 5.322 \times 10^{-13} \)}, \textcolor{black}{\( a_2 = 2.575 \times 10^{-05} \)}, and \( k = 7.12 \times 10^{3} \). The system has initial conditions: \( z(0) = 0 \), \( \dot{z}(0) = 0 \), \( \phi(0) = 0 \), \( \dot{\phi}(0) = 0 \), and \( \hat{u}(0) = 0 \). The simulation results demonstrate that the proposed ESC successfully achieves source (height) seeking, steering the flapping system towards $z_{\text{desired}}=1$. After reaching the desired point, the ESC provides viabrational stabilization for hovering with average control about 0 as observed and discussed in \cite{enatural_hovering2024}. Note that $\bm{f}(z, \dot{z}, \phi, \dot{\phi})$ in \eqref{eq:flap} is quadratic almost everywhere, which may suggest smoothing analysis similar to \cite{enatural_hovering2024} to calculate a smoothed approximation of the VOC averaged system, which we did not perform here. \textcolor{black}{We also note that the flapping application we solve here cannot be solved using the ESC approach on mechanical systems \cite{suttner2022extremum} (see also \eqref{eq:litrature_mechanical_system}) since we have to have one perturbation signal (wing-flapping) and since aerodynamic forces in flapping flight are quadratic in velocities.}

\section{Conclusion}
In this paper, we presented a novel ESC approach that provides a model-free, stable, real-time controlled vibrational stabilization for a class of mechanical systems. Inspired by flapping insects mechanics, our ESC is operable with \textit{only} one perturbation signal \textcolor{black}{which can be more natural to many mechanical systems}. \textcolor{black}{Our ESC also admits generalized forces that are quadratic in velocities, which is broader than the literature on ESC in mechanical systems, and in line with the literature on vibrational stabilization. We also were able to show the applicability of our ESC approach for generalized forces that are cubic in velocities, suggesting \textcolor{black}{even broader applications and} potential extensions in the future.} We successfully tested our ESC on multiple mechanical systems \textcolor{black}{with comparisons against the literature} and included a novel application of 1D source (height) seeking by a flapping system.

\bibliography{references} 

\begin{thebibliography}{10}
\providecommand{\url}[1]{#1}
\csname url@rmstyle\endcsname
\providecommand{\newblock}{\relax}
\providecommand{\bibinfo}[2]{#2}
\providecommand\BIBentrySTDinterwordspacing{\spaceskip=0pt\relax}
\providecommand\BIBentryALTinterwordstretchfactor{4}
\providecommand\BIBentryALTinterwordspacing{\spaceskip=\fontdimen2\font plus
\BIBentryALTinterwordstretchfactor\fontdimen3\font minus \fontdimen4\font\relax}
\providecommand\BIBforeignlanguage[2]{{%
\expandafter\ifx\csname l@#1\endcsname\relax
\typeout{** WARNING: IEEEtran.bst: No hyphenation pattern has been}%
\typeout{** loaded for the language `#1'. Using the pattern for}%
\typeout{** the default language instead.}%
\else
\language=\csname l@#1\endcsname
\fi
#2}}

\bibitem{ariyur2003real}
K.~B. Ariyur and M.~Krstic, \emph{Real-time optimization by extremum-seeking control}.\hskip 1em plus 0.5em minus 0.4em\relax John Wiley \& Sons, 2003.

\bibitem{KRSTICMain}
M.~Krstić and H.-H. Wang, ``Stability of extremum seeking feedback for general nonlinear dynamic systems,'' \emph{Automatica}, vol.~36, no.~4, pp. 595--601, 2000.

\bibitem{guay2015time}
M.~Guay and D.~Dochain, ``A time-varying extremum-seeking control approach,'' \emph{Automatica}, vol.~51, pp. 356--363, 2015.

\bibitem{scheinker2024100}
A.~Scheinker, ``100 years of extremum seeking: A survey,'' \emph{Automatica}, vol. 161, p. 111481, 2024.

\bibitem{yilmaz2023exponential}
C.~T. Yilmaz, M.~Diagne, and M.~Krstic, ``Exponential extremum seeking with unbiased convergence,'' in \emph{2023 62nd IEEE Conference on Decision and Control (CDC)}.\hskip 1em plus 0.5em minus 0.4em\relax IEEE, 2023, pp. 6749--6754.

\bibitem{khalil2002nonlinear}
H.~K. Khalil and J.~W. Grizzle, \emph{Nonlinear systems}.\hskip 1em plus 0.5em minus 0.4em\relax Prentice hall Upper Saddle River, NJ, 2002, vol.~3.

\bibitem{Maggia2020higherOrderAvg}
M.~Maggia, S.~A. Eisa, and H.~E. Taha, ``On higher-order averaging of time-periodic systems: reconciliation of two averaging techniques,'' \emph{Nonlinear Dynamics}, vol.~99, no.~1, pp. 813--836, Jan 2020.

\bibitem{DURR2013}
H.-B. Dürr, M.~S. Stanković, C.~Ebenbauer, and K.~H. Johansson, ``Lie bracket approximation of extremum seeking systems,'' \emph{Automatica}, vol.~49, no.~6, pp. 1538--1552, 2013.

\bibitem{VectorFieldGRUSHKOVSKAYA2018}
V.~Grushkovskaya, A.~Zuyev, and C.~Ebenbauer, ``On a class of generating vector fields for the extremum seeking problem: Lie bracket approximation and stability properties,'' \emph{Automatica}, vol.~94, pp. 151--160, 2018.

\bibitem{pokhrel2023higher}
S.~Pokhrel and S.~A. Eisa, ``Higher order lie bracket approximation and averaging of control-affine systems with application to extremum seeking,'' \emph{provesionally accepted in Automatica (available as arXiv preprint arXiv:2310.07092)}, 2023.

\bibitem{blekhman2000vibrational}
I.~I. Blekhman, \emph{Vibrational mechanics: nonlinear dynamic effects, general approach, applications}.\hskip 1em plus 0.5em minus 0.4em\relax World Scientific, 2000.

\bibitem{bullo2002averaging}
F.~Bullo, ``Averaging and vibrational control of mechanical systems,'' \emph{SIAM Journal on Control and Optimization}, vol.~41, no.~2, pp. 542--562, 2002.

\bibitem{nobili2023vibrational}
A.~M. Nobili, Y.~Qin, C.~A. Avizzano, D.~S. Bassett, and F.~Pasqualetti, ``Vibrational stabilization of complex network systems,'' in \emph{2023 American Control Conference (ACC)}.\hskip 1em plus 0.5em minus 0.4em\relax IEEE, 2023, pp. 1980--1985.

\bibitem{taha2020vibrational}
H.~E. Taha, M.~Kiani, T.~L. Hedrick, and J.~S. Greeter, ``Vibrational control: A hidden stabilization mechanism in insect flight,'' \emph{Science robotics}, vol.~5, no.~46, p. eabb1502, 2020.

\bibitem{tahmasian2018averaging}
S.~Tahmasian, ``On averaging and vibrational control of mechanical systems with multifrequency inputs,'' \emph{Journal of Dynamic Systems, Measurement, and Control}, vol. 140, no.~11, p. 111007, 2018.

\bibitem{tahmasian2018averaging1}
S.~Tahmasian, D.~W. Allen, and C.~A. Woolsey, ``On averaging and input optimization of high-frequency mechanical control systems,'' \emph{Journal of Vibration and Control}, vol.~24, no.~5, pp. 937--955, 2018.

\bibitem{agravcev1979exponential}
A.~Agra{\v{c}}ev and R.~V. Gamkrelidze, ``The exponential representation of flows and the chronological calculus,'' \emph{Mathematics of the USSR-Sbornik}, vol.~35, no.~6, p. 727, 1979.

\bibitem{suttner2022extremum}
R.~Suttner, ``Extremum-seeking control for a class of mechanical systems,'' \emph{IEEE Transactions on Automatic Control}, vol.~68, no.~2, pp. 1200--1207, 2022.

\bibitem{suttner2023extremum}
R.~Suttner and M.~Krsti{\'c}, ``Extremum seeking control for fully actuated mechanical systems on lie groups in the absence of dissipation,'' \emph{Automatica}, vol. 152, p. 110945, 2023.

\bibitem{grushkovskaya2021extremum}
V.~Grushkovskaya and C.~Ebenbauer, ``Extremum seeking control of nonlinear dynamic systems using lie bracket approximations,'' \emph{International Journal of Adaptive Control and Signal Processing}, vol.~35, no.~7, pp. 1233--1255, 2021.

\bibitem{observerESC}
S.~M. Mousavi and M.~Guay, ``An observer-based extremum seeking controller design for a class of second-order nonlinear systems,'' \emph{IEEE Control Systems Letters}, 2024.

\bibitem{enatural_hovering2024}
A.~A. Elgohary and S.~A. Eisa, ``Hovering flight in flapping insects and hummingbirds: A natural real-time and stable extremum seeking feedback system,'' \emph{arXiv preprint arXiv:2402.04985}, 2024.

\bibitem{MDCLUC2025}
A.~Elgohary and S.~Eisa, ``Extremum seeking in vibrational stabilization examples,'' \url{https://github.com/MDCL-UC/Extremum-seeking-in-vibrational-stabilization-examples}, 2025, accessed: 15-Feb-2025.

\bibitem{moreau2000practical}
L.~Moreau and D.~Aeyels, ``Practical stability and stabilization,'' \emph{IEEE Transactions on Automatic Control}, vol.~45, no.~8, pp. 1554--1558, 2000.

\bibitem{eisa2023}
S.~A. Eisa and S.~Pokhrel, ``Analyzing and mimicking the optimized flight physics of soaring birds: A differential geometric control and extremum seeking system approach with real time implementation,'' \emph{SIAM Journal on Applied Mathematics}, pp. S82--S104, 2023.

\end{thebibliography}
\bibliographystyle{IEEEtran} 

\end{document}